 \documentstyle[12pt]{article}
\newsavebox{\toy}
\savebox{\toy}{\framebox[0.65em]{\rule{0cm}{1ex}}}
\newcommand{\QED}{\usebox{\toy}\end{demo}}
\newenvironment{property}%
{\begin{list}{}{\setlength{\rightmargin}{0pt}%
\setlength{\itemsep}{0pt}}}{\end{list}}
\newlength{\templength}
\newcommand{\bp}{\setlength{\templength}{\labelwidth}%
\setlength{\labelwidth}{2em}\begin{property}}
\newcommand{\ep}{\end{property}\setlength{\labelwidth}{\templength}}
\newtheorem{theorem}{\indent Theorem}[section]
\newtheorem{lemma}[theorem]{\indent Lemma}
\newtheorem{proposition}[theorem]{\indent Proposition}
\newtheorem{corollary}[theorem]{\indent Corollary}
\newtheorem{remark}{\indent Remark}[section]
\newtheorem{definition}{\indent Definition}[section]
\newcommand{\Thm}[1]{Theorem \ref{Thm.#1}}

\newcommand{\Cor}[1]{Corollary \ref{Cor.#1}}

\newcommand{\Theorem}[1]{\begin{theorem}\label{Thm.#1}}
\newcommand{\Lemma}[1]{\begin{lemma}\label{Lem.#1}}
\newcommand{\Proposition}[1]{\begin{proposition}\label{Prop.#1}}
\newcommand{\Corollary}[1]{\begin{corollary}\label{Cor.#1}}
\newcommand{\Remark}[1]{\begin{remark}\label{Rem.#1}\rm}
\newcommand{\Definition}[1]{\begin{definition}\label{Def.#1}\rm}
\newcommand{\bd}{\begin{displaymath}}
\newcommand{\ed}{\end{displaymath}}
\newcommand{\bdn}{\begin{equation}}
\newcommand{\bdnl}{\begin{equation}\label}
\newcommand{\edn}{\end{equation}}
\newcommand{\barray}{\begin{array}}
\newcommand{\earray}{\end{array}}
\newcommand{\bds}{\begin{description}}
\newcommand{\eds}{\end{description}}
\newcommand{\bcenter}{\begin{center}}
\newcommand{\ecenter}{\end{center}}
\newcommand{\bflushright}{\begin{flushright}}
\newcommand{\bflushleft}{\begin{flushleft}}
\newcommand{\eflushright}{\end{flushright}}
\newcommand{\eflushleft}{\end{flushleft}}
\newcommand{\bdnn }{\begin{eqnarray*}}
\newcommand{\ednn }{\end{eqnarray*}}
\newcommand{\bdmn}{\begin{eqnarray}}
\newcommand{\edmn}{\end{eqnarray}}

\newcommand{\nn}{\nonumber}
\newcommand{\SS}[1]{\section{#1}\setcounter{equation}{0}}

\newcounter{biblio}
\newenvironment{references}%
{\begin{list}{[\arabic{biblio}]}{\usecounter{biblio}%
\setlength{\leftmargin}{2.5em}\setlength{\rightmargin}{0pt}%
\setlength{\labelwidth}{2em}\setlength{\itemsep}{0pt}}}{\end{list}}
\newcommand{\References}%
{\vspace{2.8ex plus .3ex minus .3ex}%
\begin{center}{\bf References}\end{center}\begin{references}}

\newcommand{\ra }{\rightarrow }
\newcommand{\lra }{\longrightarrow }

\newcommand{\tl }{\widetilde }
\newcommand{\ov }{\overline }

\newcommand{\ssum}{\mbox{$\sum$}}
\newcommand{\half}{\mbox{$\frac{1}{2}$}}

\newcommand{\vvs}{\vspace{2ex}}
\newcommand{\vs}{\vspace{1ex}}

\newcommand{\lef}{\left}
\newcommand{\rig}{\right}

\newcommand{\8}{\infty}

\newcommand{\6}{\partial}

\newcommand{\sub}{\subset}

\newcommand{\bsh}{\backslash}

\newcommand{\R}{{\bf R}}
\newcommand{\Z}{{\bf Z}}

\newcommand{\zz}{{\bf Z}^2}

\renewcommand{\a}{\alpha}
\renewcommand{\b}{\beta}
\newcommand{\gm}{\gamma}
\newcommand{\Gm}{\Gamma}

\newcommand{\del}{\delta}
\newcommand{\D}{\Delta}

\newcommand{\e}{\varepsilon}
\newcommand{\ve}{\epsilon}

\newcommand{\Th}{\Theta}

\newcommand{\lm}{\lambda}
\newcommand{\Lm}{\Lambda}

\newcommand{\m}{\mu}

\newcommand{\s}{\sigma}

\renewcommand{\t}{\tau}

\newcommand{\vp}{\varphi}

\newcommand{\w}{\omega}
\newcommand{\om}{\omega}
\newcommand{\W}{\Omega}
\newcommand{\Om }{\Omega}

\newcommand{\cC }{{\cal C}}

\newcommand{\epty}{\emptyset}

\makeatletter
\def\section{\@startsection{section}{1}{\z@}{-3.5ex plus -1ex minus 
 -.2ex}{2.3ex plus .2ex}{\bf}}
\def\subsection{\@startsection{subsection}{2}{\z@}{-3.25ex plus -1ex minus 
 -.2ex}{1.5ex plus .2ex}{\bf}}
\makeatother

\newcommand{\Ham }{H_{\Lm (l)}^{\om }}
\newcommand{\gibbsl }{\m_{\Lm (l),\om}^{\beta }}
\newcommand{\Omb }{\Om_{\rm b}}
\newcommand{\Oml }{\Om_{\Lm (l)}}
\newcommand{\gap }{{\rm gap}}
\newcommand{\pri}{\prime}
\begin{document}
%  \\[30pt]

\bcenter 
\large{\bf The spectral gap of the 2-D stochastic 
Ising model with mixed boundary conditions}

\normalsize 

\vvs 

Preliminary Draft

\vvs Kenneth S. Alexander\footnote{ Department of Mathematics,
University of Southern California, Los Angeles, CA 90089-1113,
USA. email: alexandr@math.usc.edu.  Research supported by NSF grant 
DMS-9802368.} 
and Nobuo Yoshida\footnote{
Division of Mathematics, Graduate School of Science, Kyoto University,
 Kyoto 606-8502, Japan.
 email:nobuo@kusm.kyoto-u.ac.jp } \\[20pt]

\ecenter 

\begin{abstract}
We establish upper bounds for the spectral gap of the stochastic Ising model
at low temperatures in an $l \times l$ box with boundary conditions which are
not purely plus or minus; specifically, we assume the magnitude of the sum
of the boundary spins over each interval of length $l$ in the boundary is
bounded by $\delta l$, where $\delta < 1$.  We show that for any such
boundary condition, when the temperature is sufficiently low (depending on
$\delta$), the spectral gap decreases exponentially in $l$.
\end{abstract}

\SS{Introduction} 
%%%%%%%%%%%%%%%%%%
\subsection{General background and heuristics}
%%%%%%%%%%%%%%%%%%
We begin with an informal description; full definitions will be given
below.  Consider the stochastic Ising model (Glauber dynamics) in an $l
\times l$ box $\Lambda(l)$, below
the critical temperature.  At equilibrium, the typical
configuration resembles one of the two infinite-volume pure phases (plus
phase or minus phase) except very near the boundary.  That is, the
equilibrium distribution
$\mu = \mu_{\Lambda(l),\omega}^{\beta}$ (at inverse temperature $\beta$,
under boundary condition $\omega$) is roughly either the plus phase, the
minus phase or a 
distributional mixture of the two.  This equilibrium may take a long time
to be reached, if the box is large.  The rate of convergence is
described by the spectral gap, denoted $\gap(\Lambda(l),\omega,\beta)$,
which is the smallest positive eigenvalue of the negative of the generator of
the dynamics. More precisely,  for $S(\cdot)$ the associated semigroup and
$\| \cdot \|_{\mu}$ the $L^{2}(\mu)$ norm, $\gap(\Lambda(l),\omega,\beta)$
is the largest constant $\Delta$ such that
\[
  \| S(t)f - \int f \, d\mu \|_{\mu} \leq \| f - \int f \, d\mu \|_{\mu}
    e^{-\Delta t} \quad 
    \mbox{for all } f \in L^{2}(\mu) \mbox{ and } t \geq 0.
\] 

For pure boundary conditions, say all plus, at subcritical temperatures the
spectral gap is believed to be of order $l^{-2}$ \cite{FH87}.  The spectral
gap can be very sensitive to the boundary condition, however.  For
example, removing as few as $O(\log l)$ plus spins near each corner of
$\Lambda(l)$ (leaving the boudary there free, or minus) yields a gap much
smaller than $l^{-2}$, and removing
$\epsilon l$ plus spins from each corner, for some positive $\epsilon$,
yields a gap which decreases exponentially in $l$ \cite{Al00}.  These
phenomena are outgrowths of the fact that the boundary conditions are not
well mixed, the free boundary or minus spins being concentrated in short
intervals at the corners.  More mixed boundary conditions are considered in
\cite{HY97}, where it is shown that if the boundary condition $\omega$
satisfies
\bdnl{hy}
|\sum_{y \in I}\w_y |  \leq  \del l/2 
\; \; \; \mbox{for every interval 
$I$ in $\6_{\rm ex}\Lm (l)$} 
\edn 
with $\delta < 1$, then 
\bdnl{expdecay}
\gap (\Lm (l), \om, \beta ) \leq B_{\ref{expdecay}}
\exp \lef( -\b l/C_{\ref{expdecay}}\rig),
\: \: \: l=1,2,\cdots ,
\edn
where $B_{\ref{expdecay}}=B_{\ref{expdecay}}(\b)>0$ and 
$C_{\ref{expdecay}}>0$. 
Here $\6_{\rm ex} \Lambda(l)$ denotes the exterior boundary; see
(\ref{bdries}).  One can allow the
boundary spins
$\omega_{y}$ to take values in the continuum $[-1,1]$, 
with $\omega_{y} = 0$ representing
the free boundary condition at site $y$.  The condition (\ref{hy}) is
somewhat restrictive, however; for example, it does not allow the long
intervals of boundary plus spins which appear in the above-mentioned results
 from
\cite{Al00}.  In this paper we establish (\ref{expdecay}) under a ``mixed
boundary'' hypothesis much weaker than (\ref{hy}).

The importance of the geometry of boundary spin locations can be seen in
comparing the result in \cite{Al00}, giving exponential decay of the gap
when $\epsilon l$ plus spins are removed at each corner, to a result of
Martinelli
\cite{Mar94} which states that when one side of the square has all-plus
boundary condition, and the other 3 sides have free boundary, at
sufficiently low temperatures one has
\bdnl{martin}
 \exp \left( -C(\beta,\epsilon)l^{\frac{1}{2} + \epsilon} \right)
 \leq \gap (\Lm (l), \om, \beta )
 \: \: \: \: \mbox{for } \epsilon > 0, \: l=1,2,\cdots.
\edn 
In the latter case there are many fewer plus spins but the gap is much
larger, meaning the convergence to the equilibrium plus phase is much
faster.

The heuristics of the gap are rooted in the ideas of energy barriers and
traps.  From certain starting configurations, to reach a typical
equilibrium configuration, one must pass through a set of configurations
for each of which the total energy is greater than either the typical
starting or equilibrium total energies.  An {\it energy barrier} is such a
set of high-energy configurations; the {\it height} of the barrier is the
typical additional energy of the barrier configurations relative to the
starting configurations.  A {\it trap} is a set of starting configurations
 from which one cannot reach equilibrium without crossing an energy barrier. 
(We do not make formal definitions here, as we will not use these concepts
other than descriptively.)  Typically one expects the gap to be
exponentially small in the height of the energy barrier that must be
crossed, for a trap 
of which the probability is ``not too small.''  Often traps are related to
the existence of macroscopic regions of the ``wrong phase,'' that is, say,
regions of minus phase when the equilibrium is purely the plus phase.  For
example, in the ``corners-removed'' context of \cite{Al00},  a trap is
formed by the configurations in which there is an ``X'' of minus phase
connecting the four free-boundary corner regions, and the height of the
associated energy barrier is proportional to the length of the corner
regions.  In the above ``three-sides-free'' example of Martinelli, however,
say with the plus spins on the right side of the square, there is no real
energy barrier because, starting from the minus phase, a region of plus
phase can sweep leftward, maintaining an approximately vertical interface,
until it covers the full square.  

Consider now a boundary condition $\omega$ which is ``well-mixed'' in the
sense that
\bdnl{mixedbdry}
|\sum_{y \in I}\w_y |  \leq  \del |I| 
\edn
for every ``sufficiently long'' interval
$I$ in the boundary of $\Lm (l)$, with $\delta < 1$,
and suppose $\omega$ favors the plus phase
(more precisely, the magnetization at the center of the square is
nonnegative.)  If the system is started entirely in the minus phase, we
expect the region of minus phase (the ``droplet'') to pull away from the
boundary and then shrink to nothing, at which time equilibrium is essentially
reached.  When the droplet initially fills $\Lambda(l)$, the energy
associated to its surface (this surface being essentially $\6_{\rm
ex}\Lambda(l)$) is at most
$8\delta l$, by (\ref{mixedbdry}). When the droplet has pulled only
slightly away from the boundary, however, the surface energy 
becomes essentially twice the droplet boundary length (provided the
temperature is very low), hence is at least about $8l$.  Thus there is an
energy barrier; the droplet will tend to stick to the boundary, meaning the
minus phase is a trap.  Though we do not make these particular heuristics
rigorous in our proofs, they are what underlie our main result.

For fixed $\omega$ satisfying (\ref{mixedbdry}),
at higher but still subcritical temperatures, one does not expect this
phenomenon of sticking to the boundary
to occur.  This is because the surface energy (appropriately defined using
surface tension and coarse-graining) of the droplet is no longer essentially
twice its length; a diagonal interface has significantly less surface energy
than combined horizontal and vertical interfaces having the same endpoints. 
This means the droplet should be able to pull away from the boundary, first
 from the corners, without the crossing of an energy barrier.  We will not
investigate this type of behavior here.

Additional existing results at subcritical temperatures include the
following.  Thomas \cite{Tho89} proved that in general dimension $d$,  
for free boundary conditions ($\omega \equiv 0$), for sufficiently 
large $\b$, 
\bdnl{thomas}
\gap (\Lm (l), \om, \beta ) \leq B_{\ref{thomas}}
\exp \lef( -\b l^{d-1}/C_{\ref{thomas}}\rig)
\: \: \: l=1,2,\cdots ,
\edn 
where $B_{\ref{thomas}}=B_{\ref{thomas}}(\b, d)>0$ and 
$C_{\ref{thomas}}=C_{\ref{thomas}}(d)>0$. 
For $d=2$, Cesi {\it et al} \cite{CGMS96} 
prove (\ref{thomas}) with $\w \equiv 0$ for all $\b >\b_c$, where 
$\b_c$ is the inverse critical temparature.  For $d = 2$, 
in contrast with (\ref{thomas}), 
it is known that 
for $\b >\b_c$ and $\w \equiv +1$, 
\bdnl{martin2}
 \exp \left( -\vp (l) \right)
 \leq \gap (\Lm (l), \om, \beta ),
 \: \: \: \: l=1,2,\cdots, 
\edn 
with a function $\vp (l)=o(l^{\half + \epsilon})$ as $l \nearrow \8$,
for all $\epsilon > 0$.  
This result was first obtained by 
F. Martinelli \cite{Mar94, Mar97}. 
More recently, 
Y. Higuchi and J. Wang \cite{HW99} showed (\ref{martin2}) 
with $\vp (l)=C(\b )(l \ln l )^{\half}$.
Schonmann \cite[Theorem 5]{Sch94} showed that
$\gap (\Lm (l), \om )$ can shrink no faster than 
exponential of $O (l^{d-1})$; specifically, 
the spectral gap has the following 
general lower bound for all $d \geq 2$ and $\b >0$:
\bdn
\underline{q}(\b)l^{-d} \exp
\left( -4\b \sum^{d-1}_{j=0}l^{j} \right) 
\leq 
\inf_{\w \in \Omb }\gap (\Lm (l), \om, \beta ),
\: \: \: l=1,2,\cdots.
\label{sch}
\edn 
Here $\underline{q}(\b)$ is a uniform lower bound for all flip rates.

\subsection{Basic definitions}
%%%%%%%%%%%%%%%%%%%%%%%%%%%%%%%%%%%%% 

\vvs 
%%%%%%%%%%%%%%%%%%%%%%%%%%%%%%
%
%   Lattice 
%
%%%%%%%%%%%%%%%%%%%%%%%%%%%%%%%
{\it The lattice. } 
For $x=(x_{1},x_{2}) \in \zz $, we will 
use both the $l_{1}$-norm 
$\|x \|_{1}=|x_{1}| + |x_{2}|$ and the
$l_{\8 }$-norm
$\|x \|_{\8 }=\max\{ |x_{1}|,|x_{2}| \}$.
A set $\Lm \sub \zz $ is said to be 
$\l_{p }$-{\it connected} ($p=1$ or $\8 $) 
if for 
each distinct $x ,y \in \Lm $, we can 
find some $ \{ x_{0}, \cdots ,x_{n}\} \sub \Lm $
with $x_{0}=x$, $x_n=y$ and 
$\| x_j-x_{j-1}\|_{p }=1$ ($j=1,\cdots ,n $).
The interior and exterior 
boundaries of a set $\Lm \sub \zz $
 will be denoted respectively 
by 
\bdmn 
\6_{\rm in}\Lm & = & \{ x \in \Lm ; 
\: \| x-y \|_{1}=1 \mbox{ for some $y \not\in \Lm $ } \}, \nn\\
\6_{\rm ex}\Lm & = & \{ y \not\in \Lm ; 
\: \| x-y \|_{1}=1 \mbox{ for some $x \in \Lm $ } \}.  \label{bdries}
\edmn 
The number of points contained in a set 
$\Lm \sub \zz$ will be denoted by $|\Lm |$. We will use the 
notation $\Lm \sub \sub \zz$ to indicate that 
$\Lm \sub \zz$ with $|\Lm |<\8 $. 
A cube with the side-length $l$ will be denoted by
\bdnl{cube}
\Lm (l)=\lef( -l/2,l/2\rig]^d \cap  \Z^d. 
\edn 
An $\l_\8$-connected subset of $\6_{\rm ex} \Lm (l)$ 
will be called an {\it interval} of $\6_{\rm ex} \Lm (l)$. 
%%%%%%%%%%%%%%%%%%%%%%%%%%%%%%%%%%%%%%%%%%%%%%%%%
%  Configurations 
%%%%%%%%%%%%%%%%%%%%%%%%%%%%%%%%%%%%%%%%%%%%%%

\vs 
{\it The configurations and the Gibbs states. }
We define two kind  of spin configurations;  
\bdnn 
 \Om_{\Lm } & = & \{ \s =(\s_x)_{x \in \Lm }\: ;
\: \s_x=\mbox{+1 or $-1$} \},\: \: \: 
\Lm \sub \sub \zz, \\
\Omb  & = & \{ \w =(\w_y)_{y \in \zz }\: ;
\: \w_y \in [-1,1]\}. 
\ednn 
We are mainly interested in $\omega_{y} \in \{ -1,0,1 \}$, but there is 
no extra work in allowing $\w_y \in [-1,1]$.
We will refer an element $\w$ of $\Omb $ 
as a boundary condition. 
The set of all real functions on $\W_{\Lm }$ is 
denoted by $\cC_{\Lm }$.
For $\Lm \sub \sub \zz$ and 
$\w \in \Omb $, the {\it Hamiltonian }
$H^{\w }_{\Lm }:  \Om_{\Lm } \ra \R $ is defined by 
$$
 H^{\w }_{\Lm }(\s )
 =-\half \sum_{\barray{c}
 \mbox{\scriptsize $x,y \in \Lm $ } \\
 \mbox{\scriptsize $\| x-y \|_{1}=1 $} 
 \earray}\s_x\s_y
 -\sum_{ \barray{c}
 \mbox{\scriptsize $x \in \Lm $, $ y \not\in \Lm $ } \\
 \mbox{\scriptsize $\| x-y \|_{1}=1 $} 
 \earray  }\s_x\w_y.
$$  
%%%%%%%%%%%%%%%%%%%%%%%%%%%%%%%%%%%%%%%%%%%%%

\vs  A {\it finite-volume Gibbs state } on $\Lm $ 
with the boundary condition $\w \in \Omb $ 
is defined to be a probability distribution 
$\m ^{\w }_{\Lm }$ on $\W_{\Lm } $, in which the 
probability of each configuration $\s \in \W_{\Lm } $ 
is given by 
$$
\m ^{\w }_{\Lm }(\{  \s \} )=  
\frac{1}{Z^{\w}_{\Lm }}\exp-\b H^{\w }_{\Lm }(\s ),
$$
where $\b >0 $ is the {\it inverse temparature} 
and $Z^{\w}_{\Lm }$ is the normalizing constant.

%%%%%%%%%%%%%%%%%%%%%%%%%%%%%%%%%%%%%%%%%%%%%%%%%%%%%%
%
%
% Stochastic Ising Models
%
%
%%%%%%%%%%%%%%%%%%%%%%%%%%%%%%%%%%%%%%%%
\vvs
{\it Stochastic Ising Models. }
For $\Lm \sub \sub \zz $, we consider a function  
$
 q_{\Lm }: \Lm \times \W_{\Lm } \times \Omb \ra \: \: ]0,\8 [
$
which satisfies the following conditions; 
\bds
\item[(i)]
Boundedness : There exist positive constants 
$\underline{q}(\b)$ 
and 
$\ov{q}(\b)$ 
such that 
\begin{equation} \label{qbound}
\underline{q}(\b)
\leq q_{\Lm }(x, \s ,\w)
\leq \ov{q}(\b),
\end{equation}
for all $\Lm \sub \sub \zz $, and 
$(x, \s , \w) \in \Lm \times \W_{\Lm } \times \Omb $.
\item[(ii)]
the Detailed Balance Condition:
\bdn
q_{\Lm }(x,\s ,\w )\exp \left( -\b H^{\w }_{\Lm }(\s ) \right)
=q_{\Lm }(x,\s^x ,\w )\exp \left( -\b H^{\w }_{\Lm }(\s^{x} ) \right), 
\label{dbc}
\edn
for all $\Lm \sub \sub \zz$ and 
$ (x,\s ,\w ) \in \Lm \times \W_{\Lm } \times \Omb $,
where $\s^x $ is the configuration obtained from 
$\s $ by replacing $\s_x$ by $-\s_x$. 
\eds 
An example of such $q_{\Lm }(x,\s ,\w )$ is given by; 
\bdnn
q_{\Lm }(x,\s ,\w )
 & = & 
 \exp \lef( -(\b /2)
 \{H^{\w }_{\Lm }(\s^{x} )-H^{\w }_{\Lm }(\s )\} \rig) \\
 & = &
 \exp \lef( -\b \s_x\lef\{ 
 \sum_{ y \in \Lm :\| x-y \|_{1}=1}\s_y
 +\sum_{ y \not\in \Lm :\| x-y \|_{1}=1  }\w_y \rig\} \rig). 
 \ednn 
Now, fix  $\Lm \sub \sub \zz $ and  $\w \in \Omb$. 
We define a linear operator 
$A^{\w }_{\Lm }:\cC_{\Lm } \ra \cC_{\Lm }$ by 
$$
A^{\w }_{\Lm }f (\s )
=\sum_{x \in \Lm }q_{\Lm }(x, \s,\w )
\{ f( \s^{x} )-f(\s ) \},
\: \: \: f \in \cC_{\Lm }. 
$$
Thus $q_{\Lm }(x, \s,\w )$ represents the rate at which the spin at $x$
flips to the opposite spin, when the configuration is $\s$. It can easily be
seen from (\ref{dbc}) that 
\bdn
-\m_{\Lm,\w}^{\b}(fA^{\w}_{\Lm }g)
  =  
 \half \sum_{x \in \Lm }
 \sum_{\s \in \W }\m_{\Lm,\w}^{\b}(\s )
 q_{\Lm }(x, \s ,\w )\{f(\s^{x})-f(\s) \} 
 \{g(\s^{x})-g(\s) \}.
\nn \edn 
  Next, we define 
 \bdn 
 \gap (\Lm ,\w,\beta )=\inf 
 \left\{ \frac{-\m_{\Lm,\w}^{\b}(fA^{\w}_{\Lm }f)}
 {\m_{\Lm,\w}^{\b}(|f-\m_{\Lm,\w}^{\b}f|^{2})}\: ;\: 
 f \in \cC_{\Lm } \right\},
\label{varcha} 
\edn 
which is the smallest positive eigenvalue of 
$-A^{\w}_{\Lm }$. 
Considering only indicator functions in (\ref{varcha}) we obtain
\begin{equation} \label{E:indicators}
  \gap(\Lambda,\omega,\beta) 
    \leq \frac{\ov{q}(\b)}{\gibbsl (\Gm)\gibbsl (\Gm^{c})}
    \sum_{x \in \Lm (l)}
    \sum_{\s \in \Gm: \s^x \not\in \Gm }\gibbsl (\s).
\end{equation}
Thus any fixed event $\Gm$ gives an upper bound for the gap.  Roughly, to
obtain a good bound one wants to choose $\Gm$ to be a trap.
%%%%%%%%%%%%%%%%%%
\subsection{Statement of main results}
%%%%%%%%%%%%%%%%%%%%%%%%%%%%%%%%%%%%%

The following is our main result, improving on the condition 
(\ref{hy}).
%%%%%%%%%%%%%%%%%%%%%%%%%%%%%%%%%%%%%%%
\Theorem{Td<1}
%%%%%%%%%%%%%%%%
Consider a stochastic Ising model on a square $\Lambda(l)$ satisfying
(\ref{qbound}) and (\ref{dbc}).  Suppose that $0 <\del <1$ and the boundary
condition 
$\w_y \in [-1,+1]$, $y \in \6_{\rm ex}\Lm (l)$  satisfies 
\bdnl{w1}
|\sum_{y \in I}\w_y |  \leq  \del |I|
\; \; \; \mbox{for every interval 
$I \sub \6_{\rm ex}\Lm (l)$ with $|I| = l$.} 
\edn 
Then, there exists $\b_0=\b_0 ( \del ) >0$ such that 
\bdn 
\gap (\Lm (l), \om, \beta )\leq B_{\ref{d<1}}
\exp (-\b l/C_{\ref{d<1}}),
\label{d<1}
\edn
for  $\b \geq \b_0 $ and $l=1,2, \ldots$, where 
$B_{\ref{d<1}}=B_{\ref{d<1}}(\b, \del )>0$ 
and $C_{\ref{d<1}}=C_{\ref{d<1}}( \del )>0$.
%%%%%%%%%%%%%%%%%%%%%%%%%%%%%%
\end{theorem}
%%%%%%%%%%%%%%%%%%%%%%%%%%%%%%%%%%%%%%%%%%%%%%%%%%%%%
Condition (\ref{w1}) is much milder than (\ref{hy}). For example,  
(\ref{w1}) allows a boundary condition which is $+1$ for 
99 \% of the boundary with  1 \% zero on each side. 
Moreover, condition (\ref{w1}) turns out 
to be optimal in the following example. For $\del >0$,  
consider a boundary condition $\w \in \Omb $ defined by   
\bdnl{bc}
\w_x=\left\{ 
\barray{ll}  
+1 & \mbox{ if $ x_1 =[ l / 2] + 1 $ 
and $\frac{-\del l}{2} < x_2 \leq \frac{\del l}{2}$,} \\
0 & \mbox{otherwise.}
\earray 
\right.
\edn 
In this example, we see the  
transition from (\ref{expdecay}) to (\ref{martin}) 
depending on the value of $\del$. 
By \Thm{Td<1}, one sees that (\ref{expdecay}) is true for all $\del <1$. 
On the other hand, it follows from \cite[Corollary 4.1]{Mar94}
that (\ref{martin}) holds true for $\del =1 $. 

\Thm{Td<1} has the following application to random boundary conditions. 
%%%%%%%%%%%%%%%%%%%%%%%%%%%%%%%%%%%%%%%%%%%%%%%%%%%%%
%
%    \Corollary{random}
%
%%%%%%%%%%%%%%%%%%%%%%%
\Corollary{random}
Suppose that $d=2$ and that 
$ \w_y \in [-1, 1] $, $y \in \Z^2 $ 
are i.i.d. random variables with the mean 
$m \in (-1,1)$. 
Then, there exists  
$\b_0 =\b_0(m)>0 $ as follows. For $\b \geq \b_0$, 
there are constants 
$B_{\ref{cor}}=B_{\ref{cor}}(\b, m )>0$ 
and $C_{\ref{cor}}=C_{\ref{cor}}( m )>0$
such that 
almost surely;
\bdn
\gap (\Lm (l), \om, \beta )\leq B_{\ref{cor}}
\exp (-\b l/C_{\ref{cor}}) 
\: \: \: 
\mbox{for } l = 1,2,...
\label{cor}
\edn
%%%%%%%%%%%%%%%%%%%%%%%%%%
\end{corollary}
%%%%%%%%%%%%%%%%%%%%%%%%%
Proof of \Cor{random} is similar to that of 
\cite[Corollary 2.2.]{HY97} and hence is omitted. 
%%%%%%%%%%%%%%%%%%%%%%%%%%%

%%%%%%%%%%%%%%%%%%%%%%%%%%%
\SS{Preliminaries for the proof of Theorem 1.1}
%%%%%%%%%%%%%%%%%%%%%%%%%%%%
\subsection{Contours}
%%%%%%%%%%%%%%%%%%%%%%%%
The set ${\bf B}$ of all bonds in $\zz$ is defined by
$$
{\bf B}=\{ \{x,y\} \sub \zz \: ; 
\:  \|x-y \|_{1} =1 \}.
$$ 
For a set $\Lm $, we define 
\bdnn
 {\bf B}_{\Lm } & = & \{ \{ x,y\} \in {\bf B} \: ; 
\:  (x,y) \in \Lm^2 \},\\
  \6 {{\bf B}_{\Lm }} & = & \{ \{x,y\} \in {\bf B} \: ; 
\:  (x,y) \in \Lm \times \Lm^c \}, \\ 
  \overline{\bf B}_\Lm & = & {\bf B}_\Lm \cup \6 {\bf B}_\Lm.
\ednn 

The \emph{dual lattice} $(\zz)^{*}$ is $\zz$ shifted by
$(\frac{1}{2},\frac{1}{2})$; sites and bonds of this lattice are called
\emph{dual sites} and \emph{dual bonds}.  $x^{*}$ denotes $x +
(\frac{1}{2},\frac{1}{2})$.  When necessary for clarity,
bonds of
$\zz$ are called \emph{regular bonds}.  To each regular bond $b$ there
is associated a unique dual bond $b^{*}$ which is its perpendicular
bisector.   For $A \subset {\bf B}$ we write $A^{*}$ for
$\{ e^{*}: e \in A \}$.   For $\gm \subset \ov{\bf B}_{\Lm}^*$ we set
\[
  V({\gamma}) = \cup_{ e = \{ x,y \}: e^{*} \in \gamma } \ \{ x,y \}, \qquad
  V_{\rm ex}(\gm) = V(\gm) \cap \6_{\rm ex}\Lm.
\]
When convenient we view bonds and dual bonds as closed intervals in 
${\bf R}^{2}$, as when referrring to a connected set of (dual) bonds.
The number of dual bonds contained in a set 
$\gm \sub {\bf B}^* $ will be denoted by $|\gm |$. 

For $x \in {\bf R}^{2}$ let $Q(x)=\prod^{2}_{j=1}
[x_{j}-\frac{1}{2}, x_{j}+\frac{1}{2}]$, and for $\Theta \subset \zz$ let
$Q(\Theta) = \cup_{x \in \Theta} \ Q(x)$.
A {\it contour} is a finite subset $\gm \sub {\bf B}^{*}$   
which is of the form $\6 Q(\Theta)$ for some finite $\Theta \subset 
\zz$ for which both $\Theta$ and $\Theta^{c}$ are $l_{1}$-connected.
The set $\Th $ is uniquely determined by $\gm $ 
and hence is denoted by $\Th (\gm)$. 
As is well known, for each $b \in {\bf B} $ and $m =1,2,\cdots $,
\bdn
\sharp \{ \gm : 
  \mbox{$\gamma$ is a contour with $|\gm |=m $ and $\gm \ni b$ } \}
\leq 3^{m-1}.
\label{count}
\edn
If a contour $\gm $ is a subset of 
$\ov{\bf B}_{\Lm}^*$ 
for some $\Lm \sub \zz $, we say 
$\gm $ is a contour in $\Lm $. For $\s \in \Om_{\Lm} $, 
$\e =+ \mbox{ or } -$ and $\Lm \sub \zz $,
an {\it ($\epsilon$)-cluster} in $\Lm$ at $\s$ is an
$l_1$-connected component of $\{ x \in {\bf Z}^2 : \s_x = \epsilon 1 \}$.
The {\it outer boundary} of a bounded subset $A$ of ${\bf R}^2$ is the unique
connected component of $\6 A$ which is contained in the closure of the
unique unbounded component of $A^c$.
A contour 
$\gm $ is  said to be an  
($\e $)-contour in $\Lm $ at $\s $ 
if $\gm$ is the outer boundary of $Q(\Th)$ for some ($\e$)-cluster $\Th$.
A contour 
$\gm $ is  said to be a contour in $\Lm $ at $\s \in \W_{\Lm }$ 
if it is either (+)-contour in $\Lm $ at $\s $ or 
($-$)-contour in $\Lm $ at $\s $.  Note that the boundary condition does
not affect whether a given $\gm$ is an ($\e$)-contour in $\Lm$, under our
definition.
%%%%%%%%%%%%%%%%%%%%%%%%%%%%%%
%%%%%%%%%%%
\subsection{Reductions}
%%%%%%%%%%%
We may assume slightly more restrictive 
condition than (\ref{w1}) to prove \Thm{Td<1}. 

First, we may replace (\ref{w1}) with 
\bdnl{w1a}
|\sum_{y \in I}\w_y |  \leq  \del_{\ref{w1a}} |I|
\; \; \; \mbox{for every interval 
$I \sub \6_{\rm ex}\Lm (l)$ with $|I| \geq l$.} 
\edn 
This is because (\ref{w1}) for a given $\delta$ implies (\ref{w1a}) with
$\delta_{\ref{w1a}} = (1 + \delta)/2$.  
We may next strengthen (\ref{w1a}) as follows: there exists 
$0< \del_{\ref{w2}} <1$ such that 
\bdnl{w2}
|\sum_{y \in I}\w_y |  \leq  \del_{\ref{w2}} |I|
\; \; \; \mbox{for every interval 
$I \sub \6_{\rm ex}\Lm (l)$ with $|I| \geq \del_{\ref{w2}} l$.} 
\edn  
In fact, suppose that the boundary condition $\w$ 
satisfies (\ref{w1a}) for some $\del <1$.  Let $\tl{\del} < 1$ satisfy
\bdnl{tldel}
\del +
(1 + \del)(1-\tl{\del})\tl{\del}^{-1} < \tl{\del}<1.
\edn
Then (\ref{w2}) with $\del_{\ref{w2}}=\tl{\del}$ holds.
To see this, take an arbirary 
interval $I \sub \6_{\rm ex}\Lm (l)$ with $|I| \geq \tl{\del} l$.
By expanding $I$, if necessary, 
we get an interval 
$\tl{I} \supset I$ such that $|\tl{I}| \geq l$ and 
$|\tl{I}\bsh I| \leq (1-\tl{\del})l 
\leq (1-\tl{\del})\tl{\del}^{-1}|I|$. 
We then have $|\tl{I}| \leq (1 + (1-\tl{\del})\tl{\del}^{-1})|I|$, so
\bdmn
\lef| \sum_{y \in I}\w_y \rig|  
 & \leq & \lef| \sum_{y \in \tl{I}}\w_y \rig|
+(1-\tl{\del})\tl{\del}^{-1}|I|, \nn \\ 
 & \leq & \del |\tl{I}| +(1-\tl{\del})\tl{\del}^{-1}|I| \nn \\ 
 & \leq & 
\lef( \del +
(1 + \del)(1-\tl{\del})\tl{\del}^{-1} \rig)|I| \nn \\
& \leq & \tl{\del}|I| \nn.
\edmn
%%%%%%%%%%%%%%%%%%%%%%%%%%%%%%%%%%%%%%%%%%%%%%%%%%%%%
\subsection{Outline of the proof of Theorem 1.1}
%%%%%%%%%%%%%%%%%%%%%%%%
We may and will assume (\ref{w2}). 
The basic strategy to prove \Thm{Td<1} is rather standard  
\cite{Tho89,HY97}. We define an event $\Gm_l \sub \Oml$ in which 
a ``large'' contour is  present; by (\ref{E:indicators}), 
\bdmn
\gap (\Lm (l), \w,\beta ) 
& \leq & 
\frac{\ov{q}(\b)}{\gibbsl (\Gm_{l})\gibbsl (\Gm_{l}^{c})}
\sum_{x \in \Lm (l)}
\sum_{\s \in \Gm_l: \s^x \not\in \Gm_l }\gibbsl (\s).
\label{gap1}
\edmn 
To prove \Thm{Td<1}, we will show that for large $\beta$, 
$\gibbsl (\Gm_{l})$ is uniformly positive in $l$ (Lemma \ref{LmGm} below)
and
$$
\gibbsl (\Gm_{l}^{c})^{-1}
\sum_{x \in \Lm (l)}
\sum_{\s \in \Gm_l, \s^x \not\in \Gm_l }\gibbsl (\s)
$$
is exponentially small in $l$ (Lemma \ref{Lm6Gm} below). 

To define $\Gm_{l}$,
we choose $\ve_{l,\w}=\pm $ as follows;
\bdn
\ve_{l,\w}=\lef\{ \barray{ll}+ & \mbox{if $\gibbsl (\s_0) \geq 0$,} \\
- & \mbox{if $\gibbsl (\s_0) < 0$.}
\earray \rig.
\edn
We fix $\del_1$ such that $\del_{\ref{w2}} <\del_1 <1$. 
The event $\Gm_l$ is defined by 
\bdnl{Gm_l}
\Gm_l= \{ \s \in \Oml \: ;\:  
C_{l}(\s )\neq \epty \},
\edn 
where 
\bdnl{C_l(s)} 
C_{l}(\s)=\left\{ \gm \: ;\: 
\barray{l}
\mbox{$\gm$ is an
  ($\ve_{l,\w}$)-contour in $\Lm (l)$ at $\s $ such that }\\ 
\mbox{ $\gm \cap \6 Q(\Lm(l)) 
\neq \epty$ and $|\gm | \geq 2\del_1 l$}   
 \earray \right\}.
\edn 

To bound (\ref{gap1}) from above, we use 
the following two lemmas.
%%%%%%%%%%%%
\begin{lemma} \label{LmGm}
%%%%%%%%%%%%%%T
Assume (\ref{w1a}).
There exists $\b_1=\b_1 (\del_{\ref{w2}} )>0$ such that 
\bdnl{mGm}
\inf_{l \geq 1} \gibbsl (\Gm_l) 
\geq \mbox{$\frac{1}{3}$}\; \; \; 
\mbox{for $\b \geq \b_1.$}
\edn 
%%%%%%%%%%%%
\end{lemma}
%%%%%%%%%%%%%%%%
%%%%%%%%%%%%
\begin{lemma} \label{Lm6Gm}
%%%%%%%%%%%%%%
Assume (\ref{w1a}).
There exists $\b_2=\b_2 (\del_{\ref{w2}} )>0$ such that 
\bdnl{m6Gm}
\sum_{x \in \Lm (l)}
\sum_{\s \in \Gm_l: \s^x \not\in \Gm_l }\gibbsl (\s) 
\leq \gibbsl (\Gm_{l}^{c})B_{\ref{m6Gm}}\exp 
\lef( -\b l /C_{\ref{m6Gm}}\rig) 
\edn 
for  $\b \geq \b_2 $ and $l=1,2, \ldots$, where 
$B_{\ref{m6Gm}}=B_{\ref{m6Gm}}( \b, \del_{\ref{w2}} )>0$ 
and $C_{\ref{m6Gm}}=C_{\ref{m6Gm}}( \del_{\ref{w2}} )>0$.
%%%%%%%%%%%%%
\end{lemma}
%%%%%%%%%%%%%%%%
\Thm{Td<1} follows immediately by plugging
 (\ref{mGm}) and (\ref{m6Gm}) into (\ref{gap1}). In fact, 
we have for $\b \geq \max \{ \b_1, \b_2\}$ that 
\bdnl{gap2}
\gap (\Lm (l), \w, \beta ) 
\leq 
3\ov{q}(\b)
B_{\ref{m6Gm}}\exp 
\lef( -\b l /C_{\ref{m6Gm}}\rig).
\edn
%%%%%%%%%%%
$\Box$

%%%%%%%%%%%%
\SS{Proof of Lemmas \ref{LmGm} and \ref{Lm6Gm}}
%%%%%%%%%%%%
\subsection{Energy estimates for contours}
%%%%%%%%%%%%%
The proofs of Lemma \ref{LmGm} and Lemma \ref{Lm6Gm} are based on energy
estimates  for contours, which we present in this subsection. 
We have to introduce additional definitions. 
The right, left, top and bottom sides of the square 
$Q(\Lm (l))$ are denoted by $F_l^{+1}, F_l^{-1}, F_l^{+2}$ and $F_l^{-2}$,
respectively.
A set of dual bonds $\gamma \subset \ov{\bf B}_{\Lm}^{*}$ is said to be 
{\it horizontally crossing} if 
\bdnl{hcro}
A \cap F_l^1 \neq \epty 
\; \; \mbox{and} \; \; A \cap F_l^{-1} \neq \epty.
\edn 
Similarly, $\gamma$ is said to be 
{\it vertically crossing} if 
\bdnl{vcro}
A \cap F_l^2 \neq \epty
\; \; \mbox{and} \; \; A \cap F_l^{-2} \neq \epty.
\edn 
The set $A$ is said to be 
{\it crossing} if it is either horizontally crossing or 
vertically crossing. 

Suppose that $\gm_1, \ldots, \gm_p$  are contours 
in $\Lm (l)$. 
We set 
\bdnl{D_T}
\D_{\gm_1, \ldots, \gm_p} \Ham (\s)  
=  \Ham (\s)-\Ham ( T_{\gm_1} \circ \cdots \circ T_{\gm_p} \s),
\; \; \; \s \in \Oml \edn 
where we have defined a map 
$T_\gm : \Oml \ra \Oml $  for a contour $\gm$ by
\bdn
(T_\gm \s)_x =\left\{
\begin{array}{ll}
-\s_x, & \mbox{if $x \in \Th (\gm)$} \\
\s_x, & \mbox{if $x \not\in \Th (\gm)$}.
\end{array}
\right.
\label{T_gm}
\edn 
Suppose that a contour $\gm$ is non-crossing.  
Then 
$\gm \cap \lef( F_l^i\cup F_l^j\rig) =\epty $  for 
some $i,j$ with $|i|=1$ and $|j|=2$. 
Then there exists a unique connected 
component $\underline{\gm}$ of 
$\gm \bsh \6 Q(\Lm (l))$ which divides $\Lm (l)$ into 
two $l_{1}$-connected components $\widetilde{\Th}$ and 
$\Lm (l) \bsh \widetilde{\Th}$ such that 
$ \Th (\gm ) \sub \widetilde{\Th}$ and 
$F_l^i\cup F_l^j 
\sub \6 Q\lef( \Lm (l) \bsh \widetilde{\Th}\rig)$. 
We define $\ov{\gm } \sub \6 Q(\Lm (l))$ and 
the interval $I(\gm ) \sub \6_{\rm ex}\Lm
(l)$  respectively by  
\bdmn 
\ov{\gm } & = & \6 Q(\Lm (l)) \cap \6 Q(\widetilde{\Th }), \label{ovgm} \\
I(\gm ) & = & V_{\rm ex}(\ov{\gm}).
\label{I(gm)} \edmn
Note that 
\bdnl{ovsup}
\ov{\gm } \supset \gm \cap \6 Q(\Lm (l)).
\edn 
Note also that bonds in $\ov{\gm }$ are in 
one-to-one correspondence with sites in $I (\gm)$ in an obvious way.

%%%%%%%%%%%%%%%
\begin{lemma} \label{Lnc}
%%%%%%%%%%%%%%%%%%
\bds
\item[a)] 
Let $\gm$ be a a non-crossing ($\ve$)-contour
at a configuration $\s \in \Oml$.  Then
\bdmn
|\gm \bsh \6 Q(\Lm (l))| & \geq & \half |\gm| 
  +\half |\gm \bsh (\6 Q(\Lm (l)) \cup \underline{\gm })|, \label{engnc0} \\
\half \D_\gm \Ham (\s)  & \geq & \half |\gm| 
  +\half |\gm \bsh (\6 Q(\Lm (l)) \cup \underline{\gm })| 
  -\ve \sum_{y \in V_{\rm ex}(\gm) }\w_y, \label{engnc} \\
\half \D_\gm \Ham (\s)  & \geq & |\underline{\gm }|
  - \ve \sum_{y \in V_{\rm ex}(\ov{\gm})} \w_y \geq 0. \label{engb}
\edmn
\item[b)] Suppose $\{\gm_j \}^p_{j=1}$ are 
($\ve$)-contours in $\Lm (l)$ at $\s$ such 
that 
\bdnl{cclm1}
\half \D_{\gm_1, \ldots, \gm_p} \Ham (\s) \; \geq \; c_1l -c_2
\; \; \; \mbox{for some $c_i \geq 0$ ($i=1,2$). }
\edn  
Then, 
\bdnl{clm1}
\half \D_{\gm_1, \ldots, \gm_p} \Ham (\s)
  \; \geq \; \mbox{$\frac{c_1}{c_1 +8}$}\sum^p_{j=1}|\gm_j |-c_2.
\edn 
\item[c)] Let $\gm$, $\gm_1$, \ldots, $\gm_p$ be non-crossing 
($\ve$)-contours 
at a configuration $\s \in \Oml$.
Suppose that condition (\ref{w2}) is satisfied 
and that $I$ is an interval in $\6_{\rm ex}\Lm (l)$ such that 
\bdmn
\cup^p_{j=1}I(\gm_j ) & \sub & I, \label{subI}\\
\del_{\ref{w2}} l & \leq & |I| \; \leq \; \sum^p_{j=1}|\underline{\gm_j}| +c,
\label{<I<}
\edmn
where $c \geq 0$ is a constant.  
Then 
\bdnl{nc}
\half \D_{\gm_1, \ldots, \gm_p} \Ham (\s)
  \; \geq \; \e_{\ref{nc}}\max \{ l, \; \sum^p_{j=1}|\gm_j | \}-c,
\edn
where the constant $\e_{\ref{nc}} >0$ depends only on $\del_{\ref{w2}}$.
\eds 
%%%%%%%%%%%%%%%%%%%%%%
\end{lemma}
%%%%%%%%%%%%%%
%%%%%%%%%%%%%%%%%%%%%
\indent Proof of part {\bf (a)}: 
%%%%%%%%%%%%%%%%%%%%%%%%%% 
We have the following relations:
\bdmn 
\half \D_\gm \Ham (\s)
  & = & |\gm \bsh \6 Q(\Lm (l))|
  -\ve \sum_{ y \in V_{\rm ex}(\gm)}\w_y ,
  \label{deco0} \\  
|I(\gm)| & = & |\ov{\gm}| \; \leq \; |\underline{\gm}|, \label{deco3} \\
|\ov{\gm} \bsh \gm| & \leq & |\gm \bsh (\6 Q(\Lm (l)) \cup 
  \underline{\gm})|.
  \label{deco4}
\edmn  
The equality (\ref{deco0}) is obvious. 
The inequalities (\ref{deco3}) and (\ref{deco4}) can be seen  
 from geometric considerations as follows. 
We decompose $\gm \bsh \6 Q(\Lm (l))$ into connected 
components $ \{ \lm_i \}_{i \geq 0}$, where $\lm_0=\underline{\gm}$.
Let $\t_0=\ov{\gm}$ and let
$ \{ \t_i \}_{i \geq 1}$ be connected components of $\ov{\gm }\bsh \gm $. 
We can arrange the enumeration  
so that $\lm_i$ and $\t_i$ 
have common endpoints for each $i \geq 0$. From this observation 
and the fact that $\gm $ is non-crossing, it follows that 
$|\lm_i | \geq |\t_i |$ for each $i \geq 0$. In particular, 
$|\lm_0 | \geq |\t_0 |$ and 
$\sum_{i \geq 1}|\t_i| \leq \sum_{i \geq 1}|\lm_i| $ which 
prove, respectively, (\ref{deco3}) and (\ref{deco4}).
By (\ref{ovsup}) and (\ref{deco3}), we have that 
$|\underline{\gm}| \geq |\gm \cap \6 Q(\Lm (l))|$ and hence that 
\bdnn
2|\gm \bsh \6 Q(\Lm (l))| 
& = & 
2 |\underline{\gm}| +2|\gm \bsh (\6 Q(\Lm (l)) \cup \underline{\gm})| \\
& \geq & |\underline{\gm}| +|\gm \cap \6 Q(\Lm (l)) | 
+2|\gm \bsh (\6 Q(\Lm (l)) \cup \underline{\gm})| \\
& = & |\gm |+|\gm \bsh (\6 Q(\Lm (l)) \cup \underline{\gm})|.
\ednn 
This proves (\ref{engnc0}), which together with (\ref{deco0})
implies (\ref{engnc}). 

On the other hand, we have by (\ref{deco4}) that 
$|\ov{\gm}\bsh \gm| \leq |\gm \bsh \6 Q(\Lm (l))|-|\underline{\gm}|$ and 
hence that   
\bdmn 
|\sum_{y \in V_{\rm ex}(\gm) }\w_y 
  -\sum_{ y \in I(\gm) }\w_y|
   & = & 
  |
  \sum_{y \in V_{\rm ex}(\ov{\gm}\bsh \gm)  }
  \w_y| \nn  \\
& \leq & 
  |\gm \bsh \6 Q(\Lm (l))|-|\underline{\gm}|. \label{ubw}
\edmn 
The inequality 
(\ref{engb}) follows from (\ref{deco0}) and (\ref{ubw}). \\
%%%%%%%%%%%%%%%%%%%%%%%%%%%
\indent Proof of part {\bf (b)}: 
Let $\a =1/(c_1+8) $. \\
Case 1: $\a  \sum^p_{j=1}|\gm_j |\leq  l$.
In this case, we obviously have that
$$
\half \D_{\gm_1, \ldots, \gm_p} \Ham (\s)  \geq c_1\a \sum^p_{j=1}|\gm_j |-c_2.
$$ 
Case 2: $\a \sum^p_{j=1}|\gm_j | \geq  l$. In this case,
\bdmn
\half \D_{\gm_1, \ldots, \gm_p} \Ham (\s)
& \geq & \ssum^p_{j=1}
\lef( |\gm_j \bsh \6 Q(\Lm (l))|-|\gm_j \cap \6 Q(\Lm (l))| \rig)\nn \\
& = & \ssum^p_{j=1}
\lef( |\gm_j |-2|\gm_j \cap \6 Q(\Lm (l))| \rig)\nn \\
& \geq &  \ssum^p_{j=1}|\gm_j | -8l \nn \\
& \geq &   (1-8\a )\ssum^p_{j=1}|\gm_j |. \nn
\edmn 
Therefore (\ref{clm1}) follows.  \\
%%%%%%%%%%%%%%%%%%%%%%%%%%%%%%%%
\indent Proof of part {\bf (c)}: 
%%%%%%%%%%%%%%%%%%%%%%%%%%%
It is enough to prove (\ref{cclm1}) with some $c_1>0$ and $c_2=c$. 
Recall that $\del_{\ref{w2}} <\del_1 <1$.\\
Case 1: $\sum^p_{j=1}|I(\gm_j )| \leq \del_1|I|$. In this case, 
\bdmn
\half \D_{\gm_1, \ldots, \gm_p} \Ham (\s) 
& \geq & \sum^p_{j=1}\lef( |\underline{\gm_j}|-|I(\gm_j )| \rig)\nn \\
& \geq & |I|-\sum^p_{j=1}|I(\gm_j )|-c \nn \\
& \geq & (1-\del_1)|I|-c \nn \\
& \geq & (1-\del_1)\del_{\ref{w2}} l-c, \nn 
\edmn 
which implies (\ref{cclm1}) with $c_1=(1-\del_1)\del_{\ref{w2}}$. \\
Case 2: $\sum^p_{j=1}|I(\gm_j )| \geq \del_1|I|$. 
We set $A=I \bsh \cup^p_{j=1}I(\gm_j )$ so that 
$|A| \leq (1-\del_1)|I|$. 
We then have by (\ref{engb}), (\ref{<I<}), (\ref{w2}) that   
\bdmn
\half \D_{\gm_1, \ldots, \gm_p} \Ham (\s) 
& \geq & 
\sum^p_{j=1}\lef( |\underline{\gm_j}|
-\ve \sum_{y \in I(\gm_j)}\w_y \rig)\nn \\
& \geq & 
|I|-c-\ve \sum_{y \in I}\w_y -|A|  \nn \\
& \geq & |I|-c-\del_{\ref{w2}} |I|-(1-\del_1)|I| \nn \\
& \geq & (\del_1 -\del_{\ref{w2}} )\del_{\ref{w2}} l-c,\nn
\edmn 
which implies (\ref{cclm1}) 
with $c_1=(\del_1 -\del_{\ref{w2}} )\del_{\ref{w2}}$. 
%%%%%%%%%%
$\Box$
%%%%%%%%%%

%%%%%%%%%%%%%%%
\Lemma{eng}
%%%%%%%%%%%%%%%%%%
Let $\gm$ be an ($\ve$)- contour in $\Lm (l)$ at a configuration $\s$. 
\bds
\item[a)] 
If $\gm$ intersects with exactly one of the sides $F_l^j$ 
($j=\pm 1, \pm 2 $), then 
\bdnl{verhol}
\D_\gm \Ham (\s) \geq \lef\{ \barray{ll}
|\{ \mbox{horizontal bonds in $\gm $} \}| & \mbox{if $j=\pm 1$,} \\
|\{ \mbox{vertical bonds in $\gm $} \}| & \mbox{if $j=\pm 2$.}
\earray \rig.
\edn 
\item[b)]
If $\Th (\gm ) \ni 0$ and $|\gm |<2l $, then
\bdnl{engc}
\D_\gm \Ham (\s)  \geq 2|\gm |/9.
\edn 
\eds 
%%%%%%%%%%%%%%%%%%%%%%
\end{lemma}
%%%%%%%%%%%%%%
%%%%%%%%%%%%%%%%%%%%%
\indent Proof of part {\bf (a)}: 
%%%%%%%%%%%%%%%%%%%%%%%%%%%
Suppose for example that 
$\gm $ intersects only with $F^1_l$. Then, 
\bdnn
\{ \mbox{horizontal bonds in $\gm $} \} & \sub & \gm \bsh \6 Q(\Lm (l)), \\
|\gm \cap F^1_l| 
& \leq & | \{ \mbox{vertical bonds in $\gm \bsh \6 \Lm (l) $} \} |.
\ednn 
Therefore, 
\bdnn
\half \D_\gm \Ham (\s) & \geq & |\{ \mbox{horizontal bonds in $\gm $} \}|  \\
&  & +| \{ \mbox{vertical bonds in $\gm \bsh \6 Q(\Lm (l))$} \} |
-|\gm \cap F^1_l| \\
& \geq & | \{ \mbox{horizontal bonds in $\gm $} \} |. 
\ednn 
%%%%%%%%%%%%%%%%%%%%%%%%%%%%%%%%
\indent Proof of part {\bf (b)}: 
%%%%%%%%%%%%%%%%%%%%%%%%%%%
If $\gm \cap \6_{\rm ex}\Lm (l)=\epty$, then (\ref{engc}) is obvious. 
We therefore assume that 
$\gm \cap \6 Q(\Lm (l)) \neq \epty$. 
Since $0 \in \Th (\gm )$ and $|\gm | <2l $, $\gm $ must intersect with 
exactly one of $F^j_l$ ($j=\pm 1, \pm 2$). 
Therefore, we see from (\ref{verhol}) that 
$ \half \D_\gm \Ham (\s) \geq l$, which implies 
(\ref{engc}) by Lemma \ref{Lnc}. 
%%%%%%%%%%%%
$\Box$
%%%%%%%%%%%%%%%
%%%%%%%%%%%%%%%%%
\subsection{Proof of Lemma 2.1}
%%%%%%%%%%%%%%%%%%%%%%%%%%%%
We  have that 
\bdmn 
\gibbsl (\Gm_l) 
 & \geq  &  
\gibbsl (\{ \s_0=\ve_{l,\w}1\}\cap \Gm_l) \nn \\
 & = & 
\gibbsl  (\{ \s_0=\ve_{l,\w}1\})-
\gibbsl  (\{ \s_0=\ve_{l,\w}1\}
\cap \Gm_l^c )  \nn \\
 & \geq & 
\half-\gibbsl  (\{ \s_0=\ve_{l,\w}1\}\cap \Gm_l^{c}).
\label{Gmcap}\edmn 
At a configuration $\s \in 
\{ \s_0=\ve_{l,\w}1\}\cap \Gm_l^{c}$, the point 
$0$ is enclosed by 
a ($\ve_{l,\w}$)-contour $\gm $ such that 
either $\gm \cap \6 Q(\Lm(l)) =\epty$ or $|\gm | <2\del_1 l$. 
If $\gm \cap \6 Q(\Lm(l)) =\epty$, then 
$$
\D_\gm H^{\w}_{\Lm (l)}(\s ) = 2 |\gm |.
$$
If otherwise, $\gm$
satisfies the condition  for (\ref{engc}). 
We therefore have in both cases
$$
\D_\gm H^{\w}_{\Lm (l)}(\s ) \geq 2|\gm |/9.
$$ 
By the standard Peierl's argument, 
\bdmn 
\gibbsl
(\{ \s_0=\ve_{l,\w}1\}\cap \Gm_l^{c}) 
  & \leq  &  
\sum_\gm
\gibbsl 
\lef\{ 
\D_\gm H^{\tl{\w}}_{\Lm (l)}(\s ) \geq 2|\gm |/9
\rig\}  \nn \\
 & \leq  &
\sum_\gm \exp \lef( -2 \b |\gm |/9 \rig),\label{sumgm}
\edmn 
where $\sum_\gm $ stands for summation over all contours $\gm$ 
which satisfies $\Th (\gm ) \ni 0$. By using 
the counting inequality (\ref{count}), we see that 
$$
\lim_{\b \nearrow \8 }\sum_\gm \exp \lef( -2 \b |\gm |/9 \rig)=0,
$$
which, together with (\ref{Gmcap}) and (\ref{sumgm}), implies 
Lemma \ref{LmGm}. 
%%%%%%%%%%%%%%%
$\Box$
%%%%%%%%%%%%
\subsection{Proof of Lemma 2.2}
%%%%%%%%%%%%%%%%%%%%%%%%%%%%
We may assume that $\ve_{l,\om} = +$.

\indent Step 1:
%%%%%%%%%%%%%%%%%%%%
Suppose that $\s \in \Gm_l$ and $\s^x \not\in \Gm_l$ 
for some $x \in \Lm (l)$.  We consider two cases separately at first: 
$\s_{x} = 1$ and $\s_{x} = -1$.

Consider first $\s_{x} = 1$.  Let $\gm$ be the outer boundary of the
(+)-cluster at $\s$ which contains $x$.  The way the transition from 
$\s \in \Gm_l$ to $\s^x \not\in \Gm_l$ occurs is that 
the set $C_l (\s )$ contains only the one element $\gm$, and 
the flipping of $\s_x$ shortens $\gm$ or separates $\gm$ from $\6 Q(\Lm(l))$
or makes $\gm $ break into new shorter contours.  Some
of these shorter contours may include dual bonds which were not part of $\gm$
at $\s$, but rather were part of ($-$)-contours inside
$\gm$ at $\s$.  We have then
\begin{equation}
  C_l (\s )  =  \{ \gm \}, \label{C=g} 
\end{equation}
\begin{equation}
  x \mbox{ is in or adjacent to } V(\gm); \label{gnix} 
\end{equation}
in fact if either
(\ref{C=g}) or (\ref{gnix}) fails, then $C_l(\s) = C_l(\s^x)$,
contradicting our assumption that $\s \in \Gm_l$ and $\s^x \not\in \Gm_l$. 
Further, there are 
(+)-contours $\gm_1, \ldots, \gm_m $ 
and $\gm_1^\pri, \ldots, \gm_n^\pri $ 
($m \geq 0$, $n \geq 0$, $1 \leq m+n \leq 4$) at the flipped configuration
$\s^x$, and ($-$)-contours $\a_1, \ldots, \a_k$ ($0 \leq k \leq 2$) inside
$\gm$ at $\s$, such that 
\bdmn
&  & \gm_j \cap \6 Q(\Lm(l)) \neq \epty, 
\; \; \; |\gm_j | < 2\del_1 l, \; \; \; 
\mbox{for $j=1,\ldots , m$,} \label{gm_j} \\
&  & \gm_j^\pri \cap \6 Q(\Lm(l)) = \epty, 
\; \; \; \mbox{for $j=1,\ldots , n$,} \label{gmpri} \\
&  & \lef( \gm \cup \lef( \cup_{j=1}^k \a_j \rig) \rig)
\triangle \lef( \lef(
\cup_{j=1}^m \gm_j \rig) 
\cup \lef( \cup_{j=1}^n \gm^\pri_j \rig)  \rig) 
\sub \6 Q(x), \label{gmtri} \\
&  & \Th (\gm ) \triangle 
\lef( \lef( \cup_{j=1}^m\Th ( \gm_j)\rig)  
\cup \lef( \cup_{j=1}^n \Th ( \gm_j^\pri ) \rig) \cup
\lef( \cup_{j=1}^k\Th ( \a_j)\rig) \rig)
= \{ x \}, \label{Thtri}
\edmn  
where $\triangle $ stands for the symmetric difference of two sets. 
Each $\a_j, \gm_j$ and $\gm_j^\pri$ must surround at least one neighbor of
$x$.  $\gm_1, \ldots, \gm_m $ 
and $\gm_1^\pri, \ldots, \gm_n^\pri $ are precisely the (+)-contours at
$\s_x$ which include bonds of $\gm$.  Let
$$
  S_- = \ssum_{j=1}^{m} |\gm_j| + \ssum_{j=1}^{n} |\gm_j^\pri |, \qquad
  S_+ = |\gm| + \ssum_{j=1}^{k} |\a_j|.
$$
Using (\ref{gmtri}) it is easy to see that
\begin{equation} 
  S_+ \leq S_- \leq S_+ + 4. \label{Scomp}
\end{equation}
We will show that
\begin{equation}
  \D_{\gm, \a_1, \ldots, \a_k} 
  H^\w_{\Lm (l)}(\s)  \geq  \e_{\ref{cl1}} S_+
  -C_{\ref{cl1}}, \label{cl1}
\end{equation}
where $\e_{\ref{cl1}}=\e_{\ref{cl1}} (\del )>0$ 
and $C_{\ref{cl1}}=C_{\ref{cl1}} (\del )>0$,
by using (\ref{gmtri}) and studying the contours $\gm_1, \ldots, \gm_m $ 
and $\gm_1^\pri, \ldots, \gm_n^\pri $.

Now consider $\s_x = -1$.  In this case, one possibility is that the flipping
of $\s_x$ connects together two or three (+)-clusters to create a
(+)-cluster which has a shorter outer boundary than the longest of the
original (+)-clusters had.  If we let $\gm$ denote the outer boundary of the
(+)-cluster of $x$ at $\s^x$, this means there are again (+)-contours 
$\gm_1, \ldots, \gm_m $ 
and $\gm_1^\pri, \ldots, \gm_n^\pri $ 
($m \geq 0$, $n \geq 0$, $1 \leq m+n \leq 2$), and
($-$)-contours $\a_1, \ldots, \a_k$ ($0 \leq k \leq 2$) inside $\gm$,
such that (\ref{gmpri})--(\ref{Thtri}) hold, but now $\gm$ and the $\a_j$
exist at $\s^x$ and the $\gm_j$ and $\gm_j^\pri$ exist at $\s$. Further, 
$|\gm| < 2 \del_1 l$, and in place of (\ref{gm_j}), 
\bdnl{gm_j2}
  \gm_j \cap \6 Q(\Lm(l)) \neq \epty, \; \; \; \mbox{for $j=1,\ldots , m$,} 
  \; \; \; |\gm_j | \geq 2\del_1 l \; \; \; 
  \mbox{for some $j$.}
\edn
(The other
possiblity when $\s_x = -1$ is that only one (+)-cluster (call it $C$) at
$\s$ is contained in the (+)-cluster of $x$ at $\s^x$, and the flipping of
$\s_x$ to 1 shortens the boundary of $C$; this may be taken as another case
of the above with $m + n = 1$ and $0 \leq k \leq 3$.)
Here in place of (\ref{C=g}) we have
\begin{equation}
  C_l (\s )  =  \{ \gm_j: 1 \leq j \leq m, |\gm_j| \geq 2 \del_1 l \}.
\label{C=g2} 
\end{equation}
Statement (\ref{gnix}) still holds, and in place of (\ref{cl1}) we will prove
\begin{equation}
  \D_{\gm_1, \ldots, \gm_m, \gm_1^\pri, \ldots, \gm_n^\pri} 
  H^\w_{\Lm (l)}(\s)  \geq  \e_{\ref{cl12}} S_-
  -C_{\ref{cl12}}, \label{cl12}
\end{equation}

It is easy to see that for fixed $x$, both for $\s_x= 1$ and for $\s_x = -1$,
the sets $\{ \gm_1, \ldots \gm_m,$ $\gm_1^\pri, \ldots \gm_n^\pri \}$ and $\{
\gm,\a_1, \ldots, \a_k \}$ uniquely determine each other.

We now turn to the proof of (\ref{cl1}) for $\s_x = 1$ and (\ref{cl12}) for
$\s_x = -1$.   For $\s_x = 1$ we have using (\ref{Thtri}) that
\bdmn
\half \D_{\gm, \a_1, \ldots, \a_k} H^\w_{\Lm (l)}(\s ) 
  & \geq & \half \D_{\gm, \a_1, \ldots, \a_k} H^\w_{\Lm (l)}(\s ) -4 \nn \\
& \geq & \half \D_{\gm_1, \ldots, \gm_m, \gm_1^\pri, \ldots, \gm_n^\pri}
  H^\w_{\Lm (l)}(\s^x )-8 \nn \\
& = & \half \D_{\gm_1, \ldots, \gm_m} H^\w_{\Lm (l)}(\s^x )
  + \half \D_{\gm_1^\pri, \ldots, \gm_n^\pri}H^\w_{\Lm (l)}(T \s^x )
  -8 \label{DT} 
\edmn
where $T=T_{\gm_1} \circ \cdots \circ T_{\gm_m}$ (Recall
(\ref{T_gm})).
Each contour in $\{ \gm_j \}$ is non-crossing, since $|\gm_j | <2\del_1 l$. 
Therefore, we see from (\ref{engnc}) that for any $0 \leq p \leq m$,  
\bdmn 
  \half \D_{\gm_1, \ldots, \gm_m} H^\w_{\Lm (l)}(\s^x )
  & \geq & 
    \half \D_{\gm_1, \ldots, \gm_p} H^\w_{\Lm (l)}(\s^x ) 
    +\ssum_{j=p+1}^m \lef( 
    \half |\gm_j |-|\gm_j \cap \6 \Lm(l) |\rig )  \nn \\
  & \geq & \half \D_{\gm_1, \ldots, \gm_p} H^\w_{\Lm (l)}(\s^x ) \nn \\
  & \geq & 0.  \label{T} 
\edmn 
On the other hand, we have
\bdnl{T'}
  \D_{\gm_1^\pri, \ldots, \gm_n^\pri}H^\w_{\Lm (l)}(T \s^x ) = 
  2\ssum^n_{j=1}|\gm^\pri_j |,
\edn 
since $\gm_j^\pri \cap \6 \Lm(l) = \epty$. 
We have as a consequence that 
\bdnl{gmgmp} 
  \half \D_{\gm, \a_1, \ldots, \a_k} H^\w_{\Lm (l)}(\s ) \; \geq \;  
  \half \D_{\gm_1, \ldots, \gm_p} H^\w_{\Lm (l)}(\s^x ) 
  +\ssum^n_{j=1}|\gm^\pri_j | -8. 
\edn
Note also that the first term on the right-hand-side of (\ref{gmgmp}) is 
non-negative by (\ref{T}).  For $\s_x = -1$, $\gm$ is non-crossing since
$|\gm| < 2 \del_1 l$, and each $\a_j$ is inside $\gm$, so it follows from
(\ref{gmtri}) that 
\bdnl{nocr}
  \6 Q(x) \cup \bigl( \cup^m_{j=1}\gm_j \bigr) \mbox{ is non-crossing in }
  Q(\Lm(l));
\edn
in particular each $\gm_j$ is non-crossing. Therefore
(\ref{DT})--(\ref{gmgmp}) remain valid but with $\s$ and $\s^x$ interchanged;
in fact we may replace (\ref{gmgmp}) with
\bdnl{gmgmp2} 
  \half \D_{\gm_1, \ldots, \gm_m, \gm_1^\pri, \ldots, \gm_n^\pri} 
  H^\w_{\Lm(l)}(\s ) \; \geq  \;
  \half \D_{\gm_1, \ldots, \gm_p} H^\w_{\Lm (l)}(\s ) 
  +\ssum^n_{j=1}|\gm^\pri_j | . 
\edn

To bound (\ref{gmgmp}) or (\ref{gmgmp2}) from
below, we pick a number
$\del_2$ such that 
$\del_{\ref{w2}} <\del_2 <\del_1$ and consider the following three cases 
separately. \\
%%%%%%%%%%%%%%%%
\indent Case 1:
%%%%%%%%%%%%%%%%
$S_+ \geq 9l$.  Here the possible energy gain along $\6 \Lm(l)$ 
when $T$ is applied is not
enough to cancel the energy reduction in the interior.  Specifically, for
$\s_x = 1$,
\bdmn
  \half \D_{\gm_1, \ldots, \gm_p} H^\w_{\Lm (l)}(\s^x ) & \geq &
    \ssum_{j+1}^m |\gm_j \backslash \6 Q(\Lm(l)) |
    - \sum_{y \in \6 \Lm(l)} |\omega_y| \nn \\
  & \geq & \ssum_{j+1}^m |\gm_j| - 8l. \label{nocancel}
\edmn
With (\ref{Scomp}) and (\ref{gmgmp}) this shows
\bdnl{case1}
  \half \D_{\gm, \a_1, \ldots, \a_k} H^\w_{\Lm (l)}(\s ) \; \geq \;
  S_+ - 8l - 8 \; \geq \; \frac{1}{9} S_+ - 8, 
\edn
which proves (\ref{cl1}).  The same argument with (\ref{gmgmp2}) replacing
(\ref{gmgmp}) proves (\ref{cl12}) when
$\s_x = -1$.  \\
%%%%%%%%%%%%%%%%
\indent Case 2: 
%%%%%%%%%%%%%%%%%%
$\ssum^m_{j=1}|\gm_j | \leq  (\del_2 /\del_1 ) S_+$. 
Consider first $\s_x = 1$.
Here by (\ref{Scomp})
\bdnl{case2}
  \ssum^n_{j=1}|\gm^\pri_j | \; \geq \; S_+ - \ssum^m_{j=1}|\gm_j |
  \; \geq \; \lef( 1 - \frac{\del_2}{\del_1} \rig) S_+,
\edn
which, together with (\ref{gmgmp}), 
proves (\ref{cl1}) in this case. Using again (\ref{Scomp}), the same
argument with (\ref{gmgmp2}) replacing (\ref{gmgmp}) proves (\ref{cl12}) when
$\s_x = -1$. \\
%%%%%%%%%%%%%%%%
\indent Case 3: 
%%%%%%%%%%%%%%%%%%
$\ssum^m_{j=1}|\gm_j | > (\del_2 /\del_1 ) S_+ $ and $S_+ < 9l$. 
By (\ref{Scomp}), (\ref{gmgmp}) and (\ref{gmgmp2}) it is enough to prove
that 
\bdnl{enoT}
\D_{\gm_1, \ldots, \gm_p} H^\w_{\Lm (l)}(\s^x ) 
\geq \e l-c
\edn 
for some  $\e>0$, $c \geq 0$ and $1 \leq p \leq m$.  
We consider several subcases as follows. \\
%%%%%%%%%%%%%%%%
\indent Case 3.1: 
%%%%%%%%%%%%%%%%%%
$\ssum^m_{j=1}|\gm_j \cap \6 Q(\Lm(l))| \leq \del_{\ref{w2}} l$. 
Consider first $\s_x = 1$.  Since
\bdnl{del2l}
  \del_2 l \; \leq \;  \half (\del_2 /\del_1 )|\gm | 
  \; \leq \; \half (\del_2 /\del_1 ) S_+ \; \leq \; \half 
  \ssum^m_{j=1}|\gm_j|, 
\edn
we have by (\ref{T}) that 
$$ 
\half \D_{\gm_1, \ldots, \gm_m} H^\w_{\Lm (l)}(\s^x )  
 \geq  (\del_2-\del_{\ref{w2}} )l 
$$
which proves (\ref{enoT}).  For $\s_x = -1$, in place of (\ref{del2l}) we
use (\ref{gm_j2}) to obtain
$$
  \del_2 l \; \leq \;  \del_1 l \; \leq \; \half 
  \ssum^m_{j=1}|\gm_j|;
$$
otherwise the argument for (\ref{enoT}) is the same.
\\
%%%%%%%%%%%%%%%%
\indent Case 3.2: 
%%%%%%%%%%%%%%%%%%
The set $\6 Q(x) \cup \bigl( \cup^m_{j=1}\gm_j \bigr)$ is non-crossing in
$\Lm (l)$ and 
\bdnl{geqdl}
  \ssum^m_{j=1}|\gm_j \cap \6 Q(\Lm(l))| \geq \del_{\ref{w2}} l.
\edn  
In this case we
have 
$$
  \lef( \cup^m_{j=1}\gm_j \rig) \cap \lef( F_l^i\cup F_l^k\rig) =\epty 
$$  
for some $i,k$ with $|i|=1$ and $|k|=2$. 
Then there exists a connected 
subset, say $\lm$,  of 
$$
  \lef( \cup^m_{j=1}\underline{\gm_j} \rig)
  \cup \6 Q(x)
$$
which divides $\Lm (l)$ into 
two connected components $\widetilde{\Th}$ and 
$\Lm (l) \bsh \widetilde{\Th}$ such that 
$ \cup^m_{j=1}\Th (\gm_j ) \sub \widetilde{\Th}$ and 
$F_l^i\cup F_l^k
\sub \6_{\rm ex}\lef( \Lm (l) \bsh \widetilde{\Th}\rig)$. 
Note that the set $I$ defined by 
$I=\6_{\rm ex}\Lm (l)\cap \6 \widetilde{\Th }$ 
is an interval. To prove (\ref{enoT}) by applying (\ref{nc}), 
let us check (\ref{subI}) and (\ref{<I<}) with $p=m$. 
We see from the construction of $I$ that 
\bdmn
  \cup^m_{j=1}I(\gm_j ) & \sub & I, \nn \\
  |I| & \leq & |\lm | \; \leq \; \sum^m_{j=1}|\underline{\gm_j}| +4. \nn
\edmn
On the other hand, we see from (\ref{geqdl}) that 
$$
  |I| \geq \sum^m_{j=1}|\gm_j \cap \6 \Lm (l)| \geq \del_{\ref{w2}} l.
$$
We therefore have (\ref{subI}) and (\ref{<I<}) with $p=m$. \\
%%%%%%%%%%%%%%%%
\indent Case 3.3: 
%%%%%%%%%%%%%%%%%%
The set $\6 Q(x) \cup \bigl( \cup^m_{j=1}\gm_j \bigr)$ is crossing in $\Lm
(l)$.  By (\ref{nocr}) this is possible only when $\s_x = 1$.
There exist $1 \leq i<j \leq 4$ and $k \in \{ 1, 2 \}$ such that 
\bdnl{121}
  \gm_i \cap F_l^k \neq \epty \; \; \; \mbox{and} \; \; \; 
  \gm_j \cap F_l^{-k} \neq \epty .
\edn 
Let us assume (\ref{121}) with $i=1$, $j=2$, and $k=1$. Then, 
the set $\gm_1 \cup \gm_2$ cannot be vertically crossing, since 
$|\gm_1|+|\gm_2| < 4 \del_1 l$ and $\gm_1 \cup \gm_2$ is already 
horizontally crossing. Let us therefore assume that 
\bdnl{1u220}
\lef( \gm_1 \cup \gm_2  \rig) \cap F_l^{-2} =\epty. 
\edn 
We are now left with two possibilities.  \\
%%%%%%%%%%%%%%%%
\indent Case 3.3.1: 
%%%%%%%%%%%%%%%%%%
$\gm_1 \cap F_l^{2} \neq \epty $ and $\gm_2 \cap F_l^{2} \neq \epty  $.
In this case, $I=I(\gm_1) \cup I(\gm_2) \cup F_l^2$ is an interval. 
To prove (\ref{enoT}) by applying (\ref{nc}), 
we will check (\ref{subI}) and (\ref{<I<}) with $p=2$. 
We obviously have 
\bdmn
\cup^2_{j=1}I(\gm_j ) & \sub & I, \nn \\
|I| & \geq & l.\nn  
\edmn
On the other hand, it is easy to see we have the following injections:
\bdmn
  I(\gm_1 ) \cap F_l^1
    & \lra & \{ \mbox{vertical dual bonds in $\underline{\gm_1} $} \}, \nn \\
  I(\gm_2 ) \cap F_l^{-1} 
    & \lra & \{ \mbox{vertical bonds in $\underline{\gm_2}$} \}, \nn \\
  F_l^2 
    & \lra & \{ \mbox{horizontal bonds in 
    $(\underline{\gm_1} \cup \underline{\gm_2}) \cup 
    \6 Q(x)$} \}. \nn 
\edmn 
We get $|I| \leq |\underline{\gm_1}| +|\underline{\gm_2}|+4$ as a consequence. 
We therefore have (\ref{subI}) and (\ref{<I<}) with $p=2$.  \\
%%%%%%%%%%%%%%%%
\indent Case 3.3.2: 
%%%%%%%%%%%%%%%%%%
$\gm_1 \cap F_l^{2} = \epty $ or 
$\gm_2 \cap F_l^{2} = \epty  $. 
Let us assume $\gm_1 \cap F_l^{2} = \epty$, 
so that $\gm_1$ does not intersect with $F_l^{j}$, $j \neq 1$. 
The distance from $x$ to $F_l^{-1}$ is 
at most $\del_1 l$, and hence 
the distance from $x$ to $F_l^{1}$ is at least $(1-\del_1)l$.
This implies that $\gm_1$  
deviates from $F_l^{1}$ at least by distance $(1-\del_1)l -1$. 
Therefore, by (\ref{verhol}), 
\bdnn
\half \D_{\gm_1} \Ham (\s) 
& \geq & |\{ \mbox{horizontal bonds in $\gm $} \}| \\
& \geq & 2(1-\del_1)l -2, 
\ednn 
which establishes (\ref{enoT}) with $p=1$. \\
%%%%%%%%%%%%%%%%%%%
\indent Step 2:
%%%%%%%%%%%%%%%%%%%%
Suppose again that $\s \in \Gm_l$ and $\s^x \notin \Gm_l$.  If $\s_x = 1$,
then every (+)-contour outside $\gm$ at $\s$ has length at most $2 \del_1
l$, and every (+)- or ($-$)-contour inside $\gm$ does not intersect $\6
Q(\Lm(l))$.  It follows that 
\bdnl{TG}
  T_{\gm} \circ T_{\a_1} \circ \cdots \circ T_{\a_k} \s \in \Gm_l^c.
\edn
Similarly if $\s_x = -1$ then
\bdnl{TG2}
  T_{\gm_1} \circ \cdots \circ T_{\gm_m} \circ 
  T_{\gm_1^\pri} \circ \cdots \circ T_{\gm_n^\pri} \s \in \Gm_l^c.
\edn
\\
%%%%%%%%%%%%
\indent Step 3: 
%%%%%%%%%%%%%
For $x \in \Lm(l)$ and $\s \in \Gm_l$ with $\s^x \notin \Gm_l$, let
${\cal C}_+ (\s,x)$ denote the set of contours $\{ \gm,\a_1, \cdots,
\a_k \}$, defined previously, if $\s_x = 1$, and let ${\cal C}_- (\s,x) = \{
\gm_1, \cdots, \gm_m, \gm_1^\pri, \cdots, \gm_n^\pri \}$ if $\s_x = -1$.  We
have by obervations made in Step 1 that 
\bdmn 
  \lefteqn{\sum_{x \in \Lm (l)}
    \sum_{\s \in \Gm_l: \s^x \not\in \Gm_l, \s_x = 1 }\gibbsl (\s) } \
    \label{st1}  \\
  & \leq & 
    \sum_{x \in \Lm (l)}\sum_{\gm, \a_1, \cdots,\a_k }
    \gibbsl \lef\{ \s: {\cal C}_+ (\s,x)=\{ \gm, \a_1, \cdots, \a_k \}, 
    \; \D_{\gm, \a_1, \ldots, \a_k} H^\w_{\Lm (l)}(\s ) \geq
    \e_{\ref{cl1}}S_+
    -C_{\ref{cl1}} \rig\} \nn
\edmn 
where $\sum_{\gm, \a_1, \cdots,\a_k }$ stands for the summation over all
possible values of ${\cal C}_+ (\s,x)$.  Now for fixed $x$ and $n$ there are
at most
$c \cdot 3^n$ possible values of ${\cal C}_+ (\s,x)$ for which $S_+ = n$.
By the standard Peierl's argument and the obervation made in
Step 2,  we can proceed as follows, provided $\b$ is sufficiently large:
\bdmn
  \lefteqn{ \sum_{x \in \Lm (l)} \ \sum_{\gm, \a_1, \cdots,\a_k } 
    \gibbsl \biggl\{ \s : {\cal C}_+ (\s,x)=\{ \gm, \a_1, \cdots, \a_k
    \}, \; \D_{\gm, \a_1, \ldots, \a_k} H^\w_{\Lm (l)}(\s ) \geq
    \e_{\ref{cl1}}S_+
    -C_{\ref{cl1}} \biggr\} } \nn \\
  & \leq &  \sum_{x \in \Lm (l)}  \sum_{\gm, \a_1, \cdots,\a_k } 
    \exp \bigl( - \beta (\e_{\ref{cl1}}S_+ - C_{\ref{cl1}} ) \bigr) \nn \\
  & \ & \qquad \qquad \qquad \qquad \cdot
    \gibbsl \biggl\{ T_{\gm} \circ T_{\a_1} \circ \cdots \circ
    T_{\a_k} \s : {\cal C}_+ (\s,x)=\{ \gm, \a_1, \cdots, \a_k
    \} \biggr\} \nn \\
  & \leq & \sum_{x \in \Lm (l)} \sum_{n \geq 2 \del_1 l} c \cdot 3^n
    \exp \bigl( - \beta (\e_{\ref{cl1}}n - C_{\ref{cl1}} ) \bigr)
    \gibbsl(\Gm_l^c) \nn \\
  & \leq & B_{\ref{st2}}\exp \lef( -\b l/C_{\ref{st2}} \rig)
    \gibbsl(\Gm_l^c).
    \label{st2} 
\edmn 
Essentially the same argument, using ${\cal C}_- (\s,x)$ and $S_-$ in place
of ${\cal C}_+ (\s,x)$ and $S_+$, shows that  
\bdnl{sx-}
  \sum_{x \in \Lm (l)}
  \sum_{\s \in \Gm_l: \s^x \not\in \Gm_l, \s_x = -1 }\gibbsl (\s) 
  \; \leq \; B_{\ref{sx-}}\exp \lef( -\b l/C_{\ref{sx-}} \rig)
  \gibbsl(\Gm_l^c).
\edn
We conclude
(\ref{m6Gm})  from (\ref{st1}), (\ref{st2}) and (\ref{sx-}).
%%%%%%%%%%%%%%%%%%%%
$\Box$
%%%%%%%%%%%%%%%%%%%%%%%%%
\footnotesize  
%%%%%%%%%%%%%%%%

\vvs 
 {\it Acknowledgement:} 
The authors thank Y. Higuchi and N. Sugimine for helpful discussions.

%%%%%%%%%%

%\newpage 
%%%%%%%%%%%%%%%

%%%%%%%%%%%%%%%%%%%%%%%
\end{document}